\documentstyle[amsfonts,11pt]{article}

\newtheorem{lemma}{Lemma}
\newtheorem{proposition}{Proposition}

\newtheorem{theorem}{Theorem}
\newtheorem{remark}{Remark}
\def\Re{\mathop{\mathrm{Re}}}

\input epsf
\def\begfig {
\begin{figure}
\small
}
\def\endfig {
\normalsize
\end{figure}
}

\begin{document}

\title{ Embedded minimal ends asymptotic to the Helicoid}

\author{David Hoffman\thanks{Hoffman was supported by research grant 
DE-FG03-95ER25250 of
the Applied Mathematical Science subprogram of the Office of Energy
Research, U.S. Department of Energy, and
by research grant DMS-95-96201 of the National Science Foundation,
Division of Mathematical Sciences. Research at MSRI is supported in
part by NSF grant DMS-90-22140.} 
\,\,\, and\,\,\, John McCuan\thanks{McCuan was supported by a 
National Science Foundation Postdoctoral Fellowship at MSRI and 
the University of California, Berkeley.  Research at MSRI
is supported in part by NSF grant DMS-9701755.}\\ 
Mathematical Sciences Research Institute \\     
1000 Centennial Drive \\                        
Berkeley CA 94720\\
\\ [24pt]}

\maketitle

\begin{abstract}
The ends of a complete embedded minimal surface of {\em finite
total curvature} are well understood (every such end is asymptotic to a
catenoid or to a plane).  We give a similar characterization for a 
large class of ends of {\em infinite total curvature}, showing that 
each such end is asymptotic to a helicoid.  The result applies, in 
particular, to the genus one helicoid and implies that it is embedded 
outside of a compact set in ${\Bbb R}^3$.
\end{abstract}

\section*{Introduction}
In \cite{HofGen}, the genus one helicoid was constructed and strong evidence
was given that it is embedded. The surface has infinite total curvature
and contains two orthogonal, intersecting straight lines. Its single
end has Weierstrass data modeled on that of the helicoid. 
It is conformally a punctured disk, on which both $dh$ and $dg/g$ have double
poles at the puncture. Moreover, these forms have no residues.  
(See the next section for
details about the Weierstrass representation and the helicoid.)
In this paper we  study this type of end.  Our main result is the 
following.

\begin{figure}[hb]
        \centerline{{
        \epsfysize=2.5in
        \epsfxsize=2.5in
        \leavevmode\epsfbox{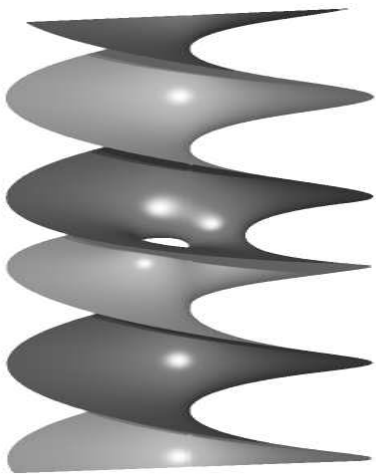}}}
\centerline{ }
\caption{The genus one helicoid.
\label{g1}}
\end{figure}

\begin{theorem} \label{mainthm} 
Let E be a complete, minimal annular end that is
conformally a punctured disk. Suppose $dg/g$ and $dh$ each have a
double pole at the puncture and that $dh$ has no residue. If E contains
a vertical ray and a horizontal ray, then a sub-end of E is embedded 
and is asymptotic to a helicoid.
\end{theorem} 
In particular, the genus one helicoid is embedded outside
of a compact set of ${\Bbb R}^3$; see Figure~\ref{g1}.  
Although our initial motivation was to establish this fact, ends with 
Weierstrass data of the type described in the Theorem~\ref{mainthm} 
are of independent interest.  
The ends of a complete embedded minimal surface of {\em finite
total curvature} are well understood (every such end is asymptotic to a
catenoid or to a plane), but a similar characterization for
embedded ends of infinite total curvature is not yet known.  A
complete minimal surface of finite total curvature is
conformally a compact Riemann surface, punctured in a finite number of
points corresponding to the ends of the minimal surface.  
The holomorphic one form $dh$ and
the meromorphic function $g$ extend meromorphically to the
compactification.  This means that $dg/g$ has, at worst, a \emph{simple
pole}, and the poles of $dh$, if there are any,  are restricted to the
punctures.  In the catenoidal case, $dh$ has a simple pole, 
which means that on a nonplanar embedded end of
finite total curvature, both $dg/g$ and $dh$ have simple poles.  If one
desires to produce an end of infinite total curvature with the conformal 
type of an isolated singularity, then $dg/g$ must have a pole of order 
two or greater. 

The end of the helicoid, on which $dg/g$ and $dh$ both have
double poles with no residue, is the only classical, embedded
annular end of infinite total curvature. It is therefore natural to ask
whether embedded, complete, minimal annular ends, on which both $dh$
and $dg/g$ have double poles, must be embedded or asymptotic to a 
helicoid.

Furthermore, as a consequence of Collin's recent positive resolution of
the generalized Nitsche Conjecture \cite{ColTop}, a complete embedded 
minimal surface
with finite topology and more than one end must have finite total
curvature. This means that a complete, embedded minimal surface with
finite topology but {\em infinite total curvature} must have just \emph{one}
end. So if we wish to restrict ourselves to looking at an embedded
annular end with infinite total curvature, which is an end of a complete
embedded minimal surface with finite topology, then it is necessarily
the only end of the surface.  Since $dh$ is holomorphic on the
punctured surface, its pole can only occur at the puncture.  
The sum of the residues of any meromorphic one-form on a compact 
Riemann surface is zero---a simple consequence of the fact that the 
union of closed curves around the punctures bounds a compact region that 
is pole-free.  In particular, 
$dh$ must have zero residue at its single pole. 
Because of this, it is reasonable to assume,
as we have in the theorem above, that the ends we study have this
property.

The authors wish to thank Harold Rosenberg, Hermann Karcher, and 
Robert Osserman for useful conversations and advice.

\section{The Weierstrass representation of the Helicoid
\label{intro-hel}} The helicoid is a simply connected, embedded, ruled
minimal surface that is generated by a family of lines, all of which
meet a fixed line orthogonally. If we assume that the fixed line is the
$x_3$-axis, then a horizontal line that sweeps out the helicoid
rotates with constant speed (as a function of $x_3$). 
From this geometric description, we will derive an analytical one:  
a Weierstrass representation of the helicoid.

A minimal surface is smooth enough so that the metric induced from 
${\Bbb R}^3$ has a conformal structure, making it a Riemann surface, $M$.  
The stereographic projection of the Gauss map, $g:M\to S^2$, is 
holomorphic, and the coordinate functions are harmonic.  
Let $h$ be the holomorphic function whose real part is $x_3$.  Even though 
$h$ is, in general, only locally defined, $dh$ is globally defined, and 
one has the Weierstrass representation 
\begin{equation}\label{general}
 X (p) = X (p_0) + \Re \int_{p_0}^p \Phi, \quad \Phi = \left ( {1\over
2} (g^{-1} -g), ~{i\over 2} ( g^{-1} +g), 1\right) \,dh. 
\end{equation}
In this representation, the integration takes place on $M$, and 
the triple $(g,dh,M)$ is referred to as the {\em Weierstrass data}.

\medskip

On the helicoid $h$ is globally defined, since the helicoid is simply 
connected.  Furthermore, from our geometric description, we know that a 
vertical helix passes through each point of the surface, so $dh$ is 
never zero.  It follows that $h$ is a diffeomorphism of the helicoid 
with $\Bbb C$.  We desire the horizontal lines of the helicoid to 
correspond to horizontal lines in $\Bbb C$.  To do this, set $z=ih$.  
Then $dh = -idz$ and, up to an additive constant, 

$$x_3(z) =  \Re \int_{z_0}^{z} dh = {\rm Im}(z), $$ 
where $z_0$ and $z$ are points in ${\Bbb C}$.

In particular, 
$X$ maps the lines ${\rm Im}(z)= {\rm constant}$ onto horizontal lines.  
Since the  vertical line meets the horizontal
lines orthogonally and the Weierstrass representation (\ref{general}) 
is conformal,
the preimage of the vertical line on the the surface is a line
${\rm Re}(z)=\alpha$ in the $z-$plane. Without loss of generality, 
we may assume
that this line is the imaginary axis, ${\rm Re}(z)=0$.

To determine $g$ as a function of $z$, note
that along the vertical axis, the normal to the helicoid is always
horizontal and rotates at a constant rate: $g(it)=e^{ic t}$, for
some nonzero real constant $c$. By scaling in ${\Bbb R}^3$ or 
by scaling $z$, we may assume, without loss of generality, that $c = 1$.
Since $g(z)$ is holomorphic, $g(z)=e^z$. 

To summarize: the Weierstrass data for the helicoid may be written 
as
$$g= e^z, ~dh = -i\,dz \hbox { on } {\Bbb C}.$$ 
We have explicitly written the helicoid as the conformal image of $\Bbb C$ 
(i.e., the Riemann sphere punctured at $\infty$).  
Up to a real translation in $z$ and a rotation in ${\Bbb R}^3$, 
this data is equivalent to 
\begin{equation}\label{eqform}
dg/g= dz, ~dh = -i\,dz \hbox { on } {\Bbb C}.  
\end{equation}  
For our purposes, we want to consider it in this form.

\section{Helicoid-type ends:  Weierstrass representation} We are
interested in understanding the behavior of a \emph{complete, minimal,
annular end, $E$, that is conformally  a punctured disk and has the
property that both $dg/g$ and $dh$ have double poles at the puncture.}
 We choose to represent the end on a neighborhood of infinity and so
may write $dh = (c_0+c_1/z+(1/z^2)f(z))dz$, where $f(z)$ is holomorphic
and bounded.  The assumption that $dh$ has a double pole at infinity
implies that $c_0\neq 0$. If ${\rm Im}(c_1) \neq 0$, then $x_3={\rm
Re}\,\int\!{dh}$ is not well defined on a neighborhood of infinity. We
will make the \emph{additional assumption} that $c_1=0$.   As was
explained in the introduction, if $E$ is an end of a complete embedded
minimal surface of finite topology, 
this is always the case. 
Under this assumption, $\zeta(z)=i\int\!dh$ is a
well defined change of variables in a neighborhood of infinity, in terms
of which 
$$
dh = -i\,d\zeta.
$$ 
Renaming the variable, we may write $dh$
and $dg/g$ in the following form:  
\begin{equation}\label{series1} 
	dg/g= (a_0+ a_1/z+(1/z^2)f(z))dz,\,\,\, ~dh = -i\,dz,
\end{equation} 
where $f(z)$ is holomorphic and
bounded on a neighborhood of infinity.  Note that the Weierstrass
representation maps any segment of ${\rm Im}(z)= {\rm constant}$ in 
the domain onto a
curve in a horizontal plane in  ${\Bbb R}^3.$

 We will prove the
following propositions.

\begin{proposition}
Let $E$ be a complete, minimal end with Weierstrass data satisfying
 (\ref{series1}). Assume that $E$ contains a vertical 
ray and a horizontal ray.
Then $E$ can be represented by Weierstrass data 
of the form
\begin{equation}\label{series2}
	dg/g= (r_0+(1/z^2)f(z))dz,\,\,\, ~dh = -i\,dz,
\end{equation} 
where 
\begin{equation}\label{form}
f(z)= \sum_{k=1}^\infty r_{2k}z^{2-2k}.
\end{equation}
The coefficients $r_{2k}$, 
$k=0,1,2,\ldots$ are real, and must 
satisfy the condition
\begin{equation}\label{coeffs}
	{\rm Res\,}_{z=0} \left(\prod_{k=0}^\infty e^{c_kz^{1-2k}} \right) =0, 
\end{equation}
 $c_k = r_{2k}/(1-2k)$.

Conversely, suppose $f$ is a bounded analytic function on a neighborhood 
of $\infty$, 
whose series representation (\ref{form}) has real coefficients 
satisfying (\ref{coeffs}). Then  the 
condition on the Weierstrass data (\ref{series2}) determines a 
minimal end containing vertical and horizontal rays, and that end 
is the unique end with these properties.   
\end{proposition}

\begin{proposition}
Any end of the type described in Proposition~1
is embedded and asymptotic to a helicoid.
\end{proposition}
\section{The proof of Proposition 1}

As noted above, the level curves of $E$ are the images under $X$ of lines 
${\rm Im}\, z = {\rm constant}$ in ${\Bbb C}$.  After a purely imaginary 
translation in $z$, we may assume without loss of generality that the 
horizontal ray is the image of a portion of the 
real axis $t\ge \rho$ (or $t\le -\rho$) for some $\rho >0$.  
On the other hand, 
the second fundamental form for $E$, as noted in \cite[(2.12)]{HofCom}, 
is given by 
$$ 
	\langle S(v)\, ,v\rangle = 
		{\rm Re}\,\left\{ {dg\over{g}}(v)\,dh(v)\right\}
$$ 
where $v$ is a tangent vector to $E$.  
In particular, along $t\ge \rho$ (or $t \le -\rho$), 
whose image is a line and hence an asymptotic curve, 
we have that $(dg/g)\,dh$ 
is purely imaginary.  Substituting from (\ref{series1}), this means that 
\begin{equation}\label{real}
	a_0 + a_1/t + (1/t^2) f(t) \in {\Bbb R}.
\end{equation}
On the other hand, $f(z)$ may be expressed as 
$$
	f(z) = \sum_{j=2}^\infty a_j z^{2-j}.
$$
In order for (\ref{real}) to hold for all large $t$, it must be the case 
that all the $a_j ,\, j\geq 0$, are real.

Through each point on the vertical ray $v$ passes a level curve of $x_3$ 
orthogonal to $v$.  Since $X$ is conformal 
and the level curves are the images of the lines 
${\rm Im}\, z = {\rm constant}$, it follows that $v$ 
is the image of a ray $\alpha + i t$ for 
$t\ge \tilde{\rho}$ (or $t\le -\tilde{\rho}$).  
By a purely real translation in $z$, we may 
assume that $\alpha =0$. 
The condition for an asymptotic curve then yields
$$
	\sum_{j=0}^\infty a_j (i\,t)^{-j} \in {\Bbb R}.
$$
Hence, $a_j= 0$ for $j$ odd. It follows from this that 
   (\ref{series2}) and 
(\ref{form}) are satisfied, and the coefficients are real. 
Because of this, we will write $r_k:=a_k$.

Next we integrate the 
expression for $dg/g$ given in (\ref{series2}) to obtain 
\begin{equation}\label{gm}
	g = e^c \prod_{k=0}^\infty e^{c_kz^{1-2k}} \equiv: A G(z),
\end{equation}
where $A=e^c$, $c$ being the constant of integration, and $c_k =
r_{2k}/(1-2k). $ 

The residue condition (\ref{coeffs}) is precisely the 
requirement that 
$$
	\int_\gamma G(z)\, dz = 0
$$
for $\gamma$ a circle centered at $0$.  

Note that $G(z)$ satisfies $G(\overline{z})=\overline{G(z)}$, and 
$G(-z) = 1/ G(z)$.
Consequently, 
\begin{equation}\label{props} 
	\overline{\int_\gamma G(z)\, dz} = \int_\gamma {1\over{G(z)}}\, dz 
					= - \int_\gamma G(z)\, dz.
\end{equation}
On the other hand, the Weierstrass representation (\ref{general}) yields 
a well defined immersion on a neighborhood of infinity provided 
\begin{equation}\label{periods}
	{\rm Re}\, \int_\gamma (g^{-1} - g)\, dh =  
			{\rm Re}\, i\int_\gamma (g^{-1} + g)\, dh = 0
\end{equation}
for a sufficiently large circle $\gamma$ centered at $0\in {\Bbb C}$.  
This is equivalent to 
$$
	\int_\gamma g^{-1}\, dh = \overline{\int_\gamma g\, dh},
$$
or 
$$
	\int_\gamma g^{-1}\, dz + \overline{\int_\gamma g\, dz} = 0
$$
since $dh = -i\, dz$.  Substituting into this identity the expression 
for $g$ from (\ref{gm}) and using (\ref{props}), we obtain 
$$
	\left(\overline{A} + {1\over{A}}\right)
		\int_\gamma G(z)dz = 0.
$$
This equation can by satisfied if and only if $\int_\gamma G(z)\, dz =0$, 
which, as noted, establishes (\ref{coeffs}).

On the other hand, Weierstrass data satisfying (\ref{series2}) and
(\ref{form}) will produce a well defined immersion (\ref{general})
provided (\ref{coeffs}) holds.  For any choice of $A$, such an
immersion always contains a horizontal ray.  To see
this, we use the formula for the planar curvature of the level
curves of a minimal surface expressed in terms of 
a Weierstrass representation with $dh= -i\, dz$ (cf. \cite{HofEmb}):
\begin{equation}\label{levelk}
	\kappa  =  {\rm
	Im}\,\left({g'\over{g}}\right)(|g|+|g|^{-1})^{-1}.  
\end{equation}
From (\ref{gm}) it is clear that $g'(t)/g(t)$ is real for $t\in {\Bbb
R}$.  Hence $\kappa \equiv 0$ on the image of any segment of the real
axis, and it follows that the image of  
$|t| \ge \rho$ is a pair of horizontal rays.

On the imaginary axis, the expression for the 
Gauss map (\ref{gm}) is unitary provided $|A| = 1$, in which case  
the image of any segment of the imaginary  
axis is a vertical line segment. 

\begin{remark}  
The lines $V$ and $H$
that contain the vertical and the horizontal rays must intersect.  If
they do not, we can produce, by Schwarz reflection about $V$, $H$, and
their images, an infinite sequence  of disjoint vertical and horizontal
lines, each of which contains a ray in the end. The images of $H$ all
lie in the same horizontal plane, forcing all their preimages in
${\Bbb C}$
to lie on a line of the form ${\rm Im} (z) = {\rm constant}$.  
By completeness, such a line
can contain the preimage of at most two disjoint rays. 

Applying rotations around the lines $V$ and $H$ above, we see that $E$
must contain both ``ends''  of $H$ and of $V$.
\end{remark}

\smallskip

Different choices of the purely imaginary constant of integration $c$
(i.e., of the unitary constant $A$) lead to rigid rotations of the same
end about the vertical ray described above.  We will show that---except
for the Weierstrass data that gives the
helicoid---if $A$ is not unitary, then any immersion defined by
(\ref{general},\ref{series2},\ref{form})  {\em contains no vertical ray}.

Suppose we have Weierstrass data in this form, defined in a
neighborhood $\cal O$ of infinity. We have shown above that the image
of ${\Bbb R}\cap \cal O$ lies on a horizontal line. The Gauss map
$g=AG$ in (\ref{gm}) is unitary on
 $i{\Bbb R}\cap {\cal O}$ if and only if $A$ is unitary there. We have
seen above that such a choice of $A$ produces an end containing  two
rays on a single vertical line, and that different unitary choices of
$A$  produce  surfaces that differ by a rotation about a vertical
axis.

What happens if we choose $|A|\neq 1$? For the end to contain vertical
rays, there must be a line $L_\alpha={\alpha+it}\subset {\Bbb C}$ such
that $g$ is unitary on $L_\alpha \cap {\cal O}$, i.e. $|G|=1/|A|$.
 From the definition of $G$ in (\ref{gm}), this is in  equivalent to
the requirement that $ c_0 z + H(z)$ has constant real part on $L_\alpha$, 
where 
$$ 
H(z)=\sum_{k=1}^\infty c_k z^{1-2k}.
$$
Since $H(z)\rightarrow 0$ as $z\rightarrow \infty$ and 
${\rm Re} \{c_0 z\} \equiv c_0 \alpha$ on $L_\alpha$, our 
requirement is in fact the geometric
condition that $H(L_\alpha \cap {\cal O})\subset i{\Bbb R}$

Inversion $z\rightarrow 1/z$ produces the function $h(z) :=
H(1/z)= \sum_{k=1}^\infty c_k z^{2k-1}$, holomorphic in a
neighborhood $\cal O'$ of  $0$, whose series expansion 
is odd with real coefficients.  In particular, $h(i{\Bbb R}\cap {\cal
O'})\subset i{\Bbb R}$ and $h(C_\alpha \cap {\cal O'})\subset i{\Bbb
R}$, where $C_\alpha$ is the circle produced by  inverting $L_\alpha$. Note
that
$C_\alpha$ is tangent to $i{\Bbb R}$ at $0$, yet they are both mapped
into $i{\Bbb R}$ by the holomorphic map $h(z)$. This is
impossible unless  $h(z)\equiv 0$, i.e. $H(z)\equiv 0$ and $g(z)
=Ae^{c_0 z}$, the Gauss map of the helicoid.  It can be put into the
standard form of \S\ref{intro-hel} by the simple change of variables
$\zeta = z +c_0^{-1} \ln|A|.$  $\Box$

\begin{remark}
The arguments in the next section rely on the
presence of a vertical ray in the surface.  It should be noted,
however, that if one chooses a  non-unitary value of $A$, then $X$
given by (\ref{general}),(\ref{series2}) is a well defined immersion
containing a horizontal ray.  It may be the case that such an end is
asymptotic to a helicoid, but the techniques of the present paper do
not apply.
\end{remark}

\begin{remark}
The condition (\ref{coeffs}) evidently
restricts the coefficients, $c_j$.  Some indication of the
nature of this restriction can be obtained by considering a simple
case.  Assume that $c_0 = 1$ and $c_k=0$ for $k>1$.  What are the
possible values of $c_1= a \ne 0$, subject to the condition 
$$
	{\rm Res}_{z=0}\left(e^{z+{a\over{z}}}\right) = 0 ?  
$$
Expanding the exponentials in series and setting the coefficient of
$1/z$ equal to zero, we obtain 
$$
	 \sum_{j=0}^\infty {(-1)^j (-a)^j\over{j!(j+1)!}} =
		{1\over{\sqrt{\xi}}}J_1(2\sqrt{-a})= 0
$$
where $J_1$ is the first Bessel function \cite[pp. 534--535]{SelCRC}.  
In particular, there is
an unbounded sequence of values, $a$, satisfying the residue condition
and yielding ends of the type described in Propositions~1~and~2.  
One member of this family of ends is illustrated in Figure~\ref{end}.  
We will return to this simple class of ends at the end of the next section.  
\begin{figure}[hb]
        \centerline{{
        \epsfysize=2.5in
        \epsfxsize=2.5in
        \leavevmode\epsfbox{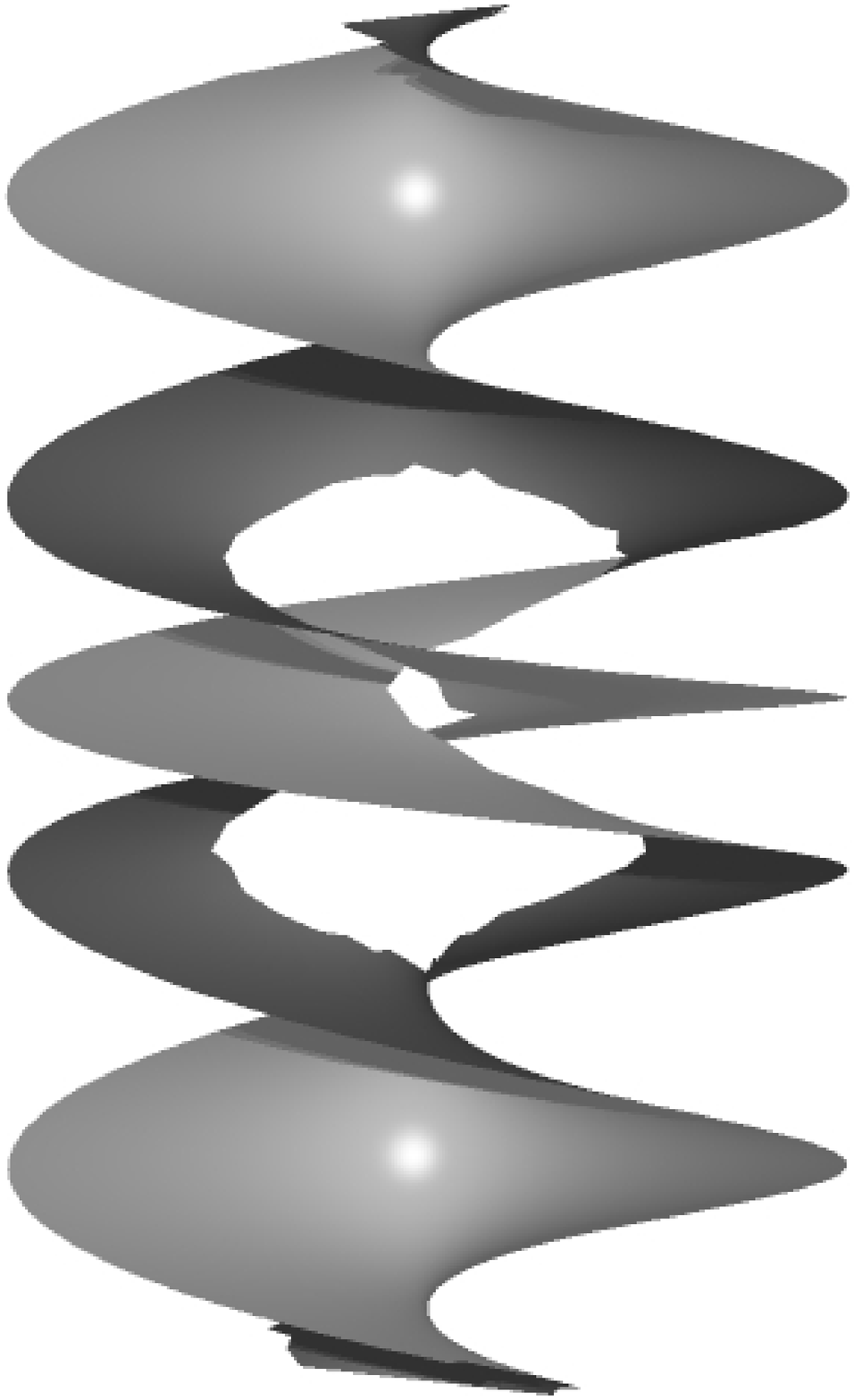}}}
\centerline{ }
\caption{The end with Gauss map $\exp(z+a/z)$ where $a$ is determined 
by the first zero of the first Bessel function.
\label{end}}
\end{figure}

\end{remark}

\section{The Proof of Proposition 2} We assume that our end $E$ is defined
on a neighborhood of $\infty$ by Weierstrass data of the form
(\ref{general}),(\ref{series2}).  Without loss of generality, we may
assume that the constant of integration $c$ is zero and---after
scaling if necessary---that $c_0=1$.  We then have
\begin{equation}\label{gm2}
	g = e^z \prod_{k=1}^\infty e^{c_kz^{1-2k}}. 
\end{equation}
 Also, we may restrict attention to
a subend of the form $X({\cal D})$ where ${\cal D} = \{ z= t+i\alpha :
|\alpha| > A \ {\rm or}\  |t|>T \}$.  We will show that if $A$ and $T$
are large enough, then the level curves $\gamma_\alpha(t)\equiv
X(t+i\alpha)$, are embedded.
  Note that for $|\alpha| > A $, $\gamma_\alpha(t)$ is defined for all
$t$, while for $|\alpha| \le A$, $\gamma_\alpha(t)$ consists of two
curves ($t>T$ or $t<-T$).

\subsection*{The End is Embedded}  We first obtain a general estimate
on the curvature of the level curves.  We may express the curvature of
$\gamma_\alpha$ in terms of the Gauss map.  From (\ref{levelk}),  we
have $$
	\kappa_\alpha(t)  =
		{\rm Im}\,\left({g'\over{g}}\right)(|g|+|g|^{-1})^{-1},
$$ valid for an immersed minimal surface with Gauss map $g$ and
$dh=-i\,dz$; the right-hand side is evaluated at $z=t+i\alpha$.  In our
case, it is convenient to use the relation $g'/g = -i(dg/g)/dh$ to
write (using (\ref{series2}) and (\ref{form})) 
\begin{eqnarray*}
	(|g|+|g|^{-1}) \kappa_\alpha(t) & = & -{\rm Re}\,\left(
					\left.{dg\over{g}}\right/
					dh\right) \\
		& = & -{\rm Re}\,\left( i\left(r_0 +
			\sum_{k=1}^\infty r_{2k}z^{-2k}\right)\right)\\
		& = & {\rm Im}\,\sum_{k=1}^\infty r_{2k}z^{-2k} .
\end{eqnarray*} 

The metric on $X({\cal D})$ is given by $ds =
(1/4)(|g|+|g|^{-1})|dz|$.  Hence, along $\gamma_\alpha(t) = X(t+i\alpha)$ 
$$
	4|\kappa_\alpha|{ds\over{dt}} \le \left|\sum_{k=1}^\infty
			r_{2k}z^{-2k}\right| \le
		\left(\sum_{k=1}^\infty{|r_{2k}|\over{|z|^{2(k-1)}}}\right)
			{1\over{|z|^2}}.
$$ 
We assume, without loss of generality, that $A$ and $T$ are large
enough so that the series $\sum r_{2k} z^{2-2k}$ is absolutely
convergent for all $z\in {\cal D}$.  All values of $\sum |r_{2k}
z^{2-2k}|$ for $z\in {\cal D}$ are clearly bounded by some $S>0$.

The total absolute curvature of $\gamma_\alpha$ is estimated by:
\begin{equation}\label{one}
	\int_{\gamma_\alpha} |\kappa_\alpha(s)|\, ds   \le 
		{S\over{4}} \int {dt\over{t^2 + \alpha^2}}\leq 
			{S\pi\over{4\alpha}}.		
\end{equation} 

In particular, when
$|\alpha|>S/4$, the total absolute curvature
 is strictly less that $\pi$. Hence, when 
 $|\alpha| \ge A_1={\rm max}\, (A,S/4)$, $\gamma_\alpha$ 
is an embedding of ${\Bbb R}$.  

\medskip

For $|\alpha|\le A_1$, we can apply the first inequality in 
(\ref{one}) to show that, for
some $T_1 >T$, $\gamma_\alpha$ is an embedding of $t>T_1$.  In fact,
$$
	\int_{T_1}^\infty |\kappa_\alpha(t)|{ds\over{dt}}\, dt \le
	{S\over{4}} \int_{T_1}^\infty {dt\over{t^2 + \alpha^2}} \le
	{S\over{4}} \int_{T_1}^\infty {dt\over{t^2}} = {S\over{4T_1}} < \pi
$$ 
for $T_1 > {\rm max}\, \{ T,{S/4\pi}\}$.  

By symmetry $\gamma_\alpha(-t)$, $t>T_1$ is also an embedded curve.  In
order to see that $\{\gamma_\alpha(t)\}$ and $\{\gamma_\alpha(-t)\}$,
$t>T_1$ are disjoint for $T_1$ large enough, we derive a general
asymptotic expression (\ref{two}) for the direction of 
$\gamma_\alpha'$, 
which will also be needed to understand the asymptotic behavior of $E$.

Since $\gamma_\alpha$ is a curve in the plane $\{x_3=\alpha\}$, there is 
a natural projection of $\gamma_\alpha$ into $\Bbb C$ by the map 
$(x_1,x_2,\alpha)\mapsto x_1+ix_2$.  It causes no confusion 
to refer to the resulting ``complex version'' of $\gamma_\alpha$ by the 
same name.  As such, the tangent vector $\gamma_\alpha'$ is a complex 
number whose argument expresses the direction of 
the tangent to the level curve $\gamma_\alpha$.  Moreover, 
${\rm arg}\,\gamma_\alpha'$ can be computed explicitly as follows:  

The horizontal projection of the Gauss map of $X({\cal D})$ at 
$z=t+i\alpha$ is a normal to the plane curve $\gamma_\alpha$ at 
$\gamma_\alpha(t)$.  Thus, $ig/|g|$ is a unit tangent vector to 
$\gamma_\alpha(t)$, whose direction is given by ${\rm arg}\,(ig)$.  In 
our case, 
\begin{equation}\label{*}
	{\rm arg}\,(g) = \alpha + {\rm Im}\, F(z),
\end{equation}
where $F(z) = \sum_{k=1}^\infty c_k z^{1-2k}$.  We may assume, without 
loss of generality, that $A$ and $T$ have been chosen so that 
$\sum |c_k z^{2-2k}|< C <\infty$ for $z\in {\cal D}$.  
In this way, 
$$ 
	|{\rm Im}\, F(z)|  \le  |z|^{-1}\sum_{k=0}^\infty |c_{k+1}| 
					|z|^{-2k} 
			\le {C\over{|z|}}.
$$ 
Combining these remarks, we obtain
\begin{equation}\label{two}
	{\rm arg}\, \gamma_\alpha'(t) = {\pi\over{2}} + \alpha + 
		{\mathcal O}\left({1\over{|z|}}\right) = 
		{\pi\over{2}} + \alpha + 
		{\mathcal O}\left(|t^2+\alpha^2|^{-1/2}\right), 
\end{equation}
where ${\mathcal O} (1/|z|) \le C/|z|$, with $C$ independent of $\alpha$ and 
$t$.

In applying (\ref{two}) it will be convenient to consider 
certain (large) neighborhoods about (the ends of) the lines 
ruling the helicoid.  Let $L_\alpha$ be the line through $0\in {\Bbb R}^2$ 
containing the direction $\vec{r}_\alpha = (\sin\alpha,\cos\alpha,0)$, and let 
$l_\alpha$ be the ray $\{ t\vec{r}_\alpha : t>0\}$.  
For a fixed $\epsilon >0$, let $N_\alpha$ be the neighborhood 
$$
	\bigcup_{|\beta - \alpha|<\epsilon} l_\beta
$$
about $l_\alpha$.  Finally, given any set $S$ in ${\Bbb R}^2$ and an 
angle $\alpha$ that will be evident from the context, we will denote by 
$\hat{S}$ the lifted set $\{ (p,\alpha)\in {\Bbb R}^3: p\in S \}$.  

\smallskip

We return now to our original use of $\gamma_\alpha$ to denote 
the level curve at height $\alpha$ in 
${\Bbb R}^3$---except where the context dictates otherwise.  
From the uniform bound in (\ref{two}) it is straightforward to show the
following:
\begin{lemma}\label{lem}
For  any $\epsilon>0$, there is some $T>0$, independent of $\alpha$, such 
that $\gamma_\alpha(t) \in \hat{N}_\alpha$ for all $t>T$.
\end{lemma}
By symmetry, $\gamma_\alpha(-t)\in \hat{N}_{\alpha+\pi}$ 
for $t$ large enough. Since $\hat{N}_{\alpha+\pi}$ and $\hat{N}_\alpha$ 
are disjoint, it is clear that the ends of each level curve are disjoint.

\subsection*{The end is asymptotic to a helicoid}    
We noted in
\S3 that, after a rigid translation of $E$ in ${\Bbb R}^3$ if
necessary,  we can assume 
that $X(i\alpha) = (0,0,\alpha)$ is the point on the standard helicoid
with Weierstrass data (\ref{eqform}). 
It is to this particular
helicoid that we will show $E$ is asymptotic.

\begin{figure}[ht]
        \centerline{{
        \epsfysize=2.5in
        \epsfxsize=2.5in
        \leavevmode\epsfbox{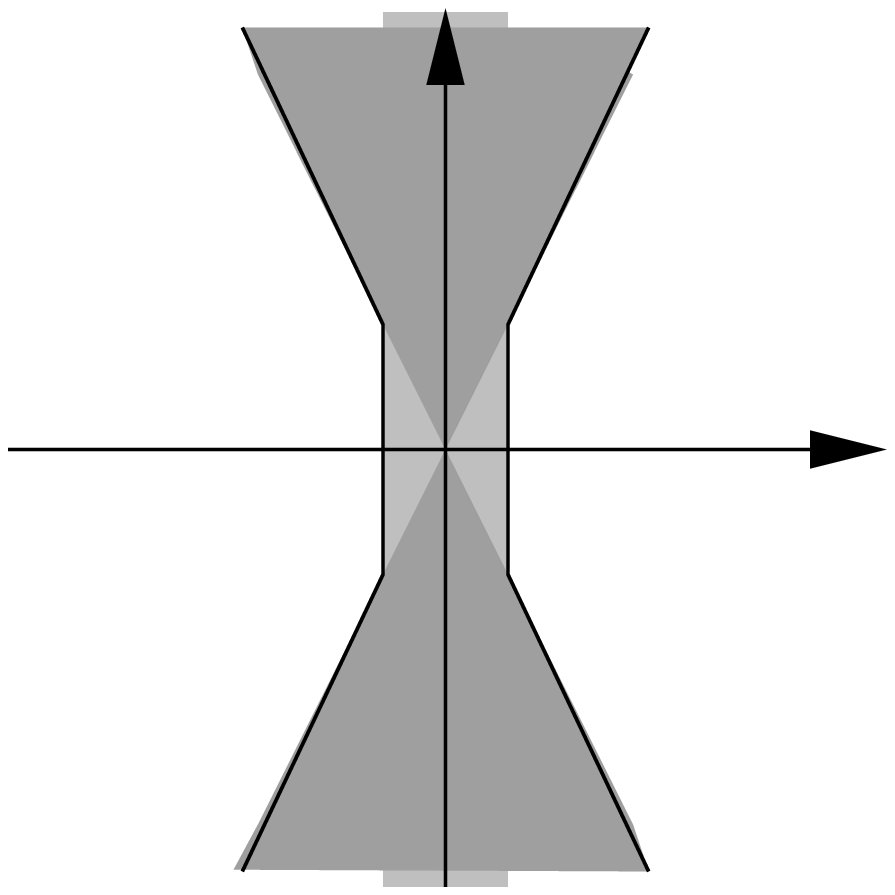}}}
\centerline{ }
\caption{A subset of $H_\epsilon \cap \{x_3 = 0 \}$.  \label{nbd}}
\end{figure}

The helicoid is ruled by lines $\hat{L_\alpha}$,  
$-\infty < \alpha <\infty$.  Let $H_\epsilon$ be the 
$\epsilon$-neighborhood of the helicoid.  It intersects the plane
$\{x_3 = 0\}$ in a set large enough to contain 
$N_0\cup N_\pi$.  See Figure~\ref{nbd}.  

Similarly, 
$$
	\hat{N}_\alpha \cup \hat{N}_{\alpha+\pi}
		\subset 
		H_\epsilon \cap \{x_3=\alpha\}.
$$
Note that a  small neighborhood of the helicoid in ${\Bbb R}^3$ contains 
quite large neighborhoods of the ruling lines in each plane 
$\{x_3=\alpha\}$.  We have already shown in Lemma~\ref{lem} that, for 
$|t|$ large enough, $X(t+ i\alpha)\in H_\epsilon$.  Let $T_2 > T_1$ be 
such that $|t| > T_2$ satisfies this condition.  We need only consider 
$\gamma_\alpha(t)$ for $|\alpha| > A_1$ and $|t| \le T_2$.  

By (\ref{two}), there is some $A_2>A_1$ such that for $|\alpha| > A_2$ 
we have 
$$
	\left| {\rm arg}\, \gamma_\alpha'(t) - \left(\alpha +\pi/2\right) 
					\right| < \epsilon.
$$
Combining this with our normalization $X(i\alpha) = \gamma_\alpha(0) = 
(0,0,\alpha)$, we see that $\gamma_\alpha(t) \in \hat{N}_\alpha\subset 
H_\epsilon$ for 
$\alpha >A_2$.  Symmetry requires, as above, that $\gamma_\alpha(-t) \in 
\hat{N}_{\alpha+\pi}$.  This establishes that the end is in $H_\epsilon$.

\begin{remark}
It cannot be expected that the level curves of 
the ends described by Theorem~1 are asymptotic to the lines ruling the 
helicoid or, in fact, to any 
lines.  If, for example, $c_1 = a$ where $J_1(2\sqrt{-a}) = 0$ as in the 
remark at the end of \S3, we see that 
$$
	\left\{ \begin{array}{ll} g(z) = e^{z+a/z} \\
				  dh = -i\, dz
		\end{array} \right.
$$
gives rise to an end $E$ with the following property:  If $\hat{E}$ is 
obtained from $E$ by a clockwise rotation about its vertical axis by an angle 
$\alpha$, i.e., $\hat{g} = e^{-i\alpha}g(z)$, then in the plane 
$\{x_3 = \alpha\}$ we have ${\rm arg}\, \hat{\gamma}_\alpha(t) \to \pi/2$, 
but up to an additive constant 
\begin{eqnarray*}
\hat{x}_1(t) & = & {\rm Re}\, \int_{i\alpha}^{i\alpha+t}{1\over{2}}\left(
			{1\over{\hat{g}}} - \hat{g}\right)\, dh \\
	& = & -{\rm Im}\, \int_0^t \sinh\left(t + 
		{{a}\over{i\alpha+t}}\right)\, dt \\
	& = & \int_0^t \sin\left( {\alpha a\over{\alpha^2+t^2}} \right)
		\cosh\left(t+{a\over{\alpha^2+t^2}} t\right)\, dt.
\end{eqnarray*}
Straightforward estimation leads to values $T,C>0$ such that for $t>T$,
$$
\int_T^t \sin\left( {\alpha a\over{\alpha^2+t^2}} \right)
		\cosh\left(t+{{a}\over{\alpha^2+t^2}} t\right)\, dt 
	\le -C \int_T^t {1\over{t^2}}e^t\, dt.
$$
Since the indefinite integral
 $\int_T^\infty t^{-2}e^t\, dt = \infty$, 
we see that $\hat{x}_1$ is unbounded.  Therefore, $\hat{\gamma}_\alpha$ 
is not asymptotic to any line, and neither is its rigid rotation 
$\gamma_\alpha$. 

\end{remark}

\bibliographystyle{plain}
\bibliography{bib/papers}

\end{document}